\documentclass[10pt]{amsart}
\usepackage{amsmath}
\usepackage{amssymb}
\usepackage{graphicx}
\usepackage{mathrsfs}

\usepackage[usenames,dvipsnames]{xcolor}

\def\beq*{\begin{eqnarray*}}
\def\eeq*{\end{eqnarray*}}

\def\be{\begin{equation}}
\def\ee{\end{equation}}

\def\a{\alpha}
\def\g{\gamma}

\newtheorem{theorem}{Theorem}

\newtheorem{remark}{Remark}

\numberwithin{subcase}{case}

\author[Pawe{\l} Hitczenko]{Pawe{\l} Hitczenko${}^\dagger$}
\address{Department of Mathematics, Drexel University, Philadelphia, 
PA  19104, USA} 
\email{phitczenko@math.drexel.edu}

\title[$\mbox{(n/log n, n/log${\bf^2}$n)}$--asymptotic normality]{A class of polynomial recurrences \\ resulting in  $\mbox{(n/log n, n/log${\bf^2}$n)}$--asymptotic \\ normality
}


\begin{document}

\begin{abstract}
We consider sequences of polynomials that satisfy differential--difference recurrences. Polynomials satisfying such recurrences frequently appear as generating polynomials of integer valued random variables that are of interest in discrete mathematics. It is, therefore, of interest to understand the properties of such polynomials and  their probabilistic consequences. We identify a class of polynomial recurrences that lead to a normal law with the expected value and the variance proportional to $n/\log n$ and $n/\log^2n$, respectively. Examples include Stirling number of the second kind and other polynomials concerning set partitions as well as polynomials related to Whitney numbers of Dowling lattices. 
\end{abstract}

\maketitle

\section{Introduction and Motivation}\label{sec:mot}
Consider a sequence of polynomials  
\[
P_n(x)=\sum_{k=0}^qp_{n,k}x^k,\quad n\ge0,
\]
where  $q=q_n$, $p_{n,k}\ge0$ and 
$\sum_{k=0}^qp_{n,k}>0$ for every $n$.  We assume $p_{n,k}=0$ for $k>q$ and often use $k\ge0$ as the range of the summation. Such polynomials are of  interest in combinatorial probability since 
\[\frac{P_n(x)}{P_n(1)}=\sum_{k\ge0}\frac{p_{n,k}}{P_n(1)}x^k\]
is the probability generating function of a non--negative,  integer valued random variable $X_n$ whose distribution function is given  by 
\begin{equation}\label{rv} {\bf P}(X_n=k)=\frac{p_{n,k}}{P_n(1)},\quad k\ge0.\end{equation}
When the underlying combinatorial objects are defined recursively, their  generating polynomials often follow recurrences of the form 
\be\label{gen-rec}P_n(x)
=a_n(x)P_{n-1}(x)+b_n(x)P_{n-1}^{'}(x)+c_n(x)P_{n-2}(x)
\ee
for specified sequences of   functions $(a_n)$, $(b_n)$, and $(c_n)$. It is therefore of interest to analyze such recurrences  and there is by now enormous literature on the subject going back to Euler, at the very least.

Many historical references and broad background are given  in  a recent paper \cite{hcd} where the authors developed the limiting theory for solutions of the above recurrence when  the polynomials are of the form:
\be\label{hk_as}a_n(x)=\a(x) n+\gamma(x),\quad b_n(x)=\beta(x)(1-x),\quad c_n(x)=0.\ee
The authors treat more  than two hundred examples when $\beta(x)\ne0$ and more than three hundred when $\beta(x)=0$ found in \cite{oeis}.
  In addition, several examples with $c_n(x)\ne0$ were discussed in Section~9.2 although in these cases the contribution of the term $c_n(x)P_{n-2}(x)$ was generally  asymptotically negligible. Within this framework the authors derived many limiting laws with various limiting distributions, including a prominently featured normal law, but also several other discrete as well as continuous distributions. The main approach was through the method of moments, but an analytic approach based on partial differential equations (PDE) and singularity analysis of the generating function  was also used. Another alternative approach is through real--rootedness of the generating polynomials (see, e.g., \cite{HL} for recent applications or \cite{hcd} for a broader discussion and more detailed references).

As was stated in \cite{hcd}, the assumption that $b_n(1)=0$  is very important for the method of moments to work.  
 Further, the assumptions made on the coefficient polynomials (particularly on $a_n(x)$) yielded a normal law with both the expected value and the variance linear in $n$, whenever the limit was Gaussian. 
 
 In the present note we concentrate on a  situation that yields a normal law with the asymptotic mean and variance proportional to $n/\log n$ and  $n/\log^2n$, respectively, under different assumptions on $a_n(x)$ and $b_n(x)$. While not nearly as extensive as the case of the Eulerian recurrences treated in \cite{hcd} it still covers a number of  cases encountered in the literature. Examples include the classical case of Stirling numbers of the second kind as well as recurrences discussed e.g. in \cite{bb}--\cite{cca}, \cite{gn}, \cite{L1}--\cite{n}, \cite{s} or \cite{wang}. 
  
 Specifically, we consider a sequence of polynomials $(P_n(x))$  satisfying the recurrence 
\be\label{rec}P_n(x)=\g(x)P_{n-1}(x)+mxP'_{n-1}(x)+(n-1)c(x)P_{n-2}(x),\quad P_0(x)=1,\ee
where $m$ is a constant and 
\be\label{gammac}\g(x)=\sum_{j=0}^k\g_jx^j,
\quad c(x)=\sum_{j=0}^\ell c_jx^j\ee
are polynomials. 
Let 
 \be\label{d}
 d:=\mbox{\rm deg}(\g(x)+c(x))\quad\mbox{and}\quad
\a_d:=[x^d]\Big( \sum_{j=1}^k\frac{\g_j}{mj}x^j+\sum_{j=1}^\ell\frac{c_j}{m^2j^2}x^j\Big).
 \ee
Then the following holds.
\begin{theorem}\label{thm}
Let $n\ge 1$. Let $X_n$ be a random variable whose probability generating function is given by \eqref{rv}
where the sequence $(P_n(x))$  satisfies the recurrence \eqref{rec}
with $m>0$. If $d$ and $\a_d$ defined in \eqref{d} satisfy $d\ge1$ and $\a_d>0$, then 
  as $n\to\infty$
\[\frac{X_n-dn/\log n}{d\sqrt n/\log n}\stackrel {\mathcal L}\longrightarrow N(0,1),\]
where \lq\lq$\stackrel {\mathcal L}\longrightarrow$\rq\rq   denotes the convergence in distribution and  $N(0,1)$ is the standard normal law.
\end{theorem}  
\begin{remark}
(i) The  conditions we imposed on the coefficient polynomials, namely, 
\[a_n(x)=\g(x), \quad b_n(x)=mx,\quad c_n(x)=(n-1)c(x),\]
seem quite restrictive. Nonetheless, in virtually all recurrences of this type $(b_n)$ do not depend on $n$. If one drops the requirement that $b_n(1)=0$, then $b_n(x)=mx$ is  common (it is also responsible for the asymptotic values of the expectation and the variance). Also, if the polynomials $c_n(x)$ are assumed to be linear functions of $n$ (with coefficients that are polynomials in $x$), assuming that they are of the form $(n-1)c(x)$ is not a serious restriction as the other term (a polynomial in $x$ not depending on $n$)  would not contribute significantly.  The assumption on $(a_n(x))$ (which corresponds to setting  $\a(x)=0$ in \eqref{hk_as}) still covers a number of cases  including  Stirling numbers of the second kind (see Section~\ref{ex} below for some specific examples). This assumption is complementary to \cite{hcd} where it was assumed $\a(1)>0$ (which, as we mentioned above, led to the limiting normal law with the expected value and the variance linear in $n$).

(ii) Other approaches to establishing asymptotic normality have been used and some of them are discussed in \cite{hcd}. One of them is based on showing that the polynomials $(P_n(x))$ have real roots only (in fact this is the case for many of the examples discussed in the references we mentioned earlier).  The real--rootedness is of interest in itself and has been studied extensively. Examples of relatively recent work in some degree of generality in that direction include \cite{DDJ,LW}. However, as pointed out in~\cite{hcd}, the real--rootedness property seems quite sensitive to variations in the coefficient polynomials $\g(x)$, $b(x)$, and $c(x)$. In addition, one still needs to show that the variance of the resulting random variables grows to infinity with $n$, which often in this context is not  a substantially easier task.  In some of the referenced papers the asymptotic normality was asserted; in others it was not. The approach requires working with specific recurrences to establish the real--rootedness. Our result provides some generality and uniformity for a class of frequently encountered polynomials. In virtually all of the examples referenced earlier, the degree of $\g(x)+c(x)$ is one resulting in the expected value and the variance asymptotic to  $n/\log n$ and $n/\log^2n$, respectively.  
\end{remark}

\section{Proof of Theorem~\ref{thm}}
Our proof proceeds along the typical lines. 
We first find the explicit form of the bivariate generating function $F(z,x)$ that encodes the probability distribution function of the underlying random variables  and then carry out the asymptotic analysis of its coefficients.

\subsection{A PDE associated with \eqref{rec}}
To derive the bivariate generating function we  consider a partial differential equation that $F(z,x)$ satisfies.  We note that in our situation the resulting PDF can be solved by the method of characteristics giving the explicit expression for $F(z,x)$.  While this approach is not new, it has not been used much in this context. For more on the method of characteristics we refer to \cite{pde} or almost any other textbook on  PDEs. One of the advantages is that it provides a simple and transparent way of deriving the bivariate generating function in cases when the resulting PDF can be solved. Our situation is particularly simple.
 Let
\be\label{Fzx}F(z,x):=\sum_{n=0}^\infty P_n(x)\frac{z^n}{n!}
.\ee
We differentiate \eqref{Fzx} with respect to $z$.  Using \eqref{rec} (and a convention that $P_j(x)=0$ whenever $j<0$) gives
\begin{align*}\frac\partial{\partial z}F(z,x)&=\sum_{n=1}^\infty
\Big\{\g(x)P_{n-1}(x)+mxP'_{n-1}(x)+(n-1)c(x)P_{n-2}(x)\Big\}\frac{z^{n-1}}{(n-1)!}\\&=
\g(x)F(z,x)+mx\frac\partial{\partial x}F(z,x)+c(x)zF(z,x)
\end{align*}
or
\be\label{pde}\frac\partial{\partial z}F(z,x)-mx\frac\partial{\partial x}F(z,x)=(\g(x)+c(x)z)F(z,x).
\ee
With $F(0,x)=P_0(x)$ (which usually is equal to 1) this is easily solved by the method of characteristics. Namely, 
by setting 
\be\label{char}\frac{dx}{dz}=-mx\ee
 PDE \eqref{pde} is reduced to the ordinary differential equation:
\[
\frac d{dz}F(z,x(z))=(\g(x(z))+zc(x(z)))F(z,x(z))
\]
whose solution is 
\[F(z,x(z))=\exp\left\{\int (\g(x(z))+zc(x(z)))dz\right\},\]
where by \eqref{char}  
$x(z)=\xi e^{- mz},
$ 
and  the parameter $\xi$ is treated as a constant of integration.  
Eliminating the  parameter $\xi=xe^{mz}$
and choosing the constant of integration so that $F(0,x)=P_0(x)$ 
gives the explicit expression for $F(z,x)$. 

In our case, using \eqref{gammac} 
we obtain
\begin{align*}&\int\Big(\g(\xi e^{-mz})+zc(\xi e^{-mz})\Big)dz=\g_0z+c_0\frac{z^2}2+C_0\\&\quad -\sum_{j=1}^k\frac{\g_j}{jm}\xi^j(e^{-jmz}+C_j)-\sum_{j=1}^\ell\frac{c_j}{jm}\xi^j\left(ze^{-jmz}+\frac{e^{-jmz}}{jm}+K_j\right),
\end{align*}
where $C_j$ and $K_j$ are integration constants. Setting $F(0,x)=1$ means that $C_0=0$, $C_j=-1$, $1\le j\le k$ and $K_j=-1/(jm)$, $1\le j\le \ell$. Thus, using $\xi=xe^{mz}$ we obtain
 \be\label{fzx}
F(z,x)=
\exp\left\{
\g_0z+c_0\frac{z^2}2+\sum_{j=1}^k\frac{\g_j}{jm}x^j(e^{jmz}-1)+\sum_{j=1}^\ell\frac{c_j}{jm}x^j\left(
\frac{e^{jmz}-1}{jm}-z\right)
\right\}
.\ee
\subsection{Asymptotics of $[z^n]F(z,x)$}
Functions given by \eqref{fzx} are amenable to the perturbation of the saddle point asymptotics. We refer to Section~IX.8 of \cite{fs}) for a discussion of the perturbation aspects in the bivariate case and to Chapter~VIII  for a detailed presentation  of the saddle point estimation.  In our case, we can invoke the general principles developed in Hayman's work
 \cite{hay} on admissibility.
  (Essentially, a function is called admissible if it is subject to the saddle point asymptotics; we refer the reader to \cite[Section~VIII.5]{fs} or the original work of Hayman \cite{hay} for more details.)  

Let $\Omega$ be a small neighborhood of  $x=1$  (in particular $x>0$ for $x\in\Omega$). We fix $x\in\Omega$ for a moment and  write the  exponent in \eqref{fzx} as 
\be\label{f-rep}f(z,x)=Q_1(z,x)+Q_2(xe^{mz})
\ee
where 
\[Q_1(z,x)=
-\sum_{j=1}^k\frac{\g_jx^j}{jm}-\sum_{j=1}^\ell\frac{c_jx^j}{j^2m^2}+\Big(\g_0-\sum_{j=1}^\ell\frac{c_jx^j}{jm}\Big)z+c_0\frac{z^
2}2
\]
and 
\[Q_2(z)=\sum_{j=1}^k\frac{\g_j}{jm}z^{j}+\sum_{j=1}^\ell\frac{c_j}{j^2m^2}z^{j}.
\]
Clearly, the function $e^{mz}$ is admissible if $m>0$ (it also follows from \cite[Theorem~X]{hay}). Since the leading coefficient of $Q_2$ is positive and $x>0$,  $Q_2(xe^{mz})$ is admissible and so is 
$Q_1(z,x)+Q_2(xe^{mz})$ by \cite[Corollary to Theorem~IX]{hay}. Thus, $f(z,x)$ is admissible and so is $F(z,x)$ by \cite[Theorem~VI]{hay}.

We will apply Hayman's result to $e^{f(z,x)}$. Choose 
$r$ 
 so that 
\be\label{sad_pt_eq}(zf_z(z,x))_{z=r}=n.
\ee
Then
\[[z^n]F(z,x)=
\frac{r^{-n}}{2\pi}\int_{-\pi}^\pi e^{f(re^{i\theta},x)-in\theta}d\theta
\sim \frac{r^{-n}e^{f(r,x)}}{\sqrt{2\pi b(r,x)}},
\]
where $b(r,x)=(zf_z(z,x)+z^2f_{zz}(z,x))_{z=r}$ and where $f_z(z,x)$, $f_{zz}(z,x)$ (and later $f_x(z,x)$, $f_{zx}(z,x)$, etc.) denote the  derivative(s) of $f$ with respect to the indicated variable(s) (not just the partial derivative(s) with respect to the first or second argument of $f$). 

By examining the argument above and noting that the dependence on $x$ is polynomial, it  is clear that the above estimates are uniform for $x\in\Omega$.  We thus infer that the  probability generating functions of random variables  encoded by $F(z,x)$ are given asymptotically by
 \[p_n(x)=\frac{F(\rho(x),x)}{F(\rho(1),1)}\left(\frac{\rho(1)}{\rho(x)}\right)^n(1+o(1)),
  \]
 where the $o(1)$ error is uniform over  $\Omega$  and $\rho(x)=\rho_n(x)$ is a positive solution of the saddle point equation \eqref{sad_pt_eq}
 for $x\in\Omega$.  In particular,  $\rho(x)$ satisfies
 \be\label{sad_pt}\rho(x)f_z(z,x)_{z=\rho(x)}=n.
 \ee
Taking the logarithms and recalling by \eqref{d} that  the leading coefficient of $Q_2(z)$ is $\a_d$ we see that
 \[\log\rho(x)+\log\left(md\a_dx^de^{md\rho(x)}\Big(1+O\Big(\frac1{xe^{m\rho(x)}}\Big)\Big)\right)=\log n.
 \]
 It follows that 
 \[md\rho(x)=\log n-\log\rho(x)-d\log x-\log md\a_d+O\Big(\frac1{xe^{m\rho(x)}}\Big).
 \]
 Successive iterations starting with $x=1$ give 
 \be\label{rho_as}\rho(x)=\frac{\log n}{md}+O(\log\log n)\ee
 with the uniform behavior in $x$ over $\Omega$.
 
 Let us denote
   \[h_n(x):=f(\rho(x),x)-n\log\rho(x),\]
so that 
\[p_n(x)=\exp(h_n(x)-h_n(1))(1+o(1)).\]
 Then, by \cite[Theorem~IX.13]{fs} on generalized quasi--powers,  the corresponding random variables are asymptotically normal  provided 
 \be\label{h'''}\frac{h'''(x)}{(h'_n(1)+h''_n(1))^{3/2}}\rightarrow0,
 \ee
 uniformly over  $\Omega$.
 Moreover, the mean and the variance are asymptotic to $h'_n(1)$ and  to $h_n'(1)+h_n''(1)$, respectively. 
 Differentiating $h_n(x)$ we get
 \[h_n'(x)=\rho'(x)f_z(z,x)_{z=\rho(x)}+f_x(z,x)_{z=\rho(x)}-n\frac{\rho'(x)}{\rho(x)}.
\]
In view of \eqref{sad_pt} 
\[\rho'(x)f_z(z,x)_{z=\rho(x)}=\frac{\rho'(x)}{\rho(x)}\rho(x)f_z(z,x)_{z=\rho(x)}=n\frac{\rho'(x)}{\rho(x)}.
\]
Thus, $h'_n(x)$  simplifies to 
 \[h_n'(x)=f_x(z,x)_{z=\rho(x)}.
 \]
 Differentiating again yields
\[
h''_n(x)=\rho'(x)f_{zx}(z,x)_{z=\rho(x)}+f_{xx}(z,x)_{z=\rho(x)}.
 \]
 Finally, implicit differentiation of \eqref{sad_pt} gives
 \[\rho'(x)f_z(z,x)_{z=\rho(x)}+\rho(x)\Big(\rho'(x)f_{zz}(z,x)_{z=\rho(x)}+f_{zx}(z,x)_{z=\rho(x)}\Big)=0.
  \]
 After rearranging the terms and putting $z=\rho(x)$
 \[\rho'(x)=\frac{-\rho(x)f_{zx}(\rho(x),x)}{f_z(\rho(x),x)+\rho(x)f_{zz}(\rho(x),x)}.
  \]
 Evaluating at $x=1$ (and writing $\rho=\rho(1)$, $h'_n=h'_n(1)$, etc.) we arrive at 
  \be\label{rhoprime}\rho'=\frac{-\rho f_{zx}(\rho,1)}{f_z(\rho,1)+\rho f_{zz}(\rho,1)}.
   \ee
Since $\rho\to\infty$ as $n\to\infty$ and $f(z,x)$ is a polynomial in $x$, $z$ and $e^z$, it is clear that the asymptotic behavior of the relevant expressions depends on the coefficients of highest power of $e^\rho$. Specifically, 
\begin{align*}f_x(\rho,1)&\sim d\a_de^{\rho md}\\
f_z(\rho,1)&\sim md\a_de^{\rho md}\\
 f_{zz}(\rho,1)&\sim m^2d^2\a_de^{\rho md}\\
 f_{zx}(\rho,1)&\sim md^2\a_de^{\rho md}\\
f_{xx}(\rho,1)&\sim d(d-1)\a_de^{\rho md}.
\end{align*}
  Further, \eqref{sad_pt} implies that
 \[\rho md\a_de^{\rho md}\sim n.\]
 Thus, using \eqref{rho_as}, 
 \be\label{hnprime}h_n'= f_x(\rho,1)\sim d\a_de^{\rho md}\sim \frac n{m\rho}\sim\frac{dn}{\log n}.
 \ee
 Similarly,
 \be\label{hndblprime}
h''_n=\rho'f_{zx}(\rho,1)+f_{xx}(\rho,1)\sim \rho'md^2\a_de^{\rho md}+d(d-1)\a_de^{\rho md}
\ee
so that
 \[h_n'+h_n''\sim (1+m\rho')d^2\a_de^{\rho md}.\]
By \eqref{rhoprime}
\[\rho'\sim\frac{-\rho md^2\a_de^{\rho md}}{md\a_de^{\rho md}+\rho m^2d^2\a_de^{\rho md}}\sim
\frac{-\rho d}{1+\rho md}
.\]
 Hence,
\begin{align*}h_n'+h_n''&\sim \left(1-\frac{\rho md}{1+m\rho d}\right)d^2\a_de^{\rho md}\sim 
\frac{1}{1+\log n}\frac dmmd\a_de^{\rho md}\\&\sim 
\frac{1}{1+\log n}\frac dm\frac n\rho \sim \frac{d^2n}{\log^2 n}
.\end{align*}
  Finally, it is clear from \eqref{hnprime}, \eqref{hndblprime} and the form of  $f(z,x)$ that \eqref{h'''} holds uniformly over a small neighborhood $\Omega$ of $x=1$. 
 This completes the proof.

 \section
 {Examples and further comments}\label{ex}
 \subsection{Set partitions of type $B_n$} Wang \cite{wang} established the normal limit law for the number of non--zero blocks in the colored set partitions of type $B_n$ (we refer to \cite {wang} for the definitions and background). This amounted to analyzing polynomial recurrences of the form
 \be\label{wang_poly}T_n(x)=(x+c)T_{n-1}(x)+mxT'_{n-1}(x),\quad n\ge 1, \quad T_0(x)=1,\ee
 where $c$ and $m$ are positive integers. 
 In order to do it, he showed that each $T_n(x)$ has real, negative roots. This implies that the resulting $X_n$ is a sum of independent indicators. He then used the formulas
 \[{\bf E} X_n =\frac{T_{n+1}(1)}{mT_n(1)}-\frac{1+c}m,\quad {\bf var}(X_n)=\frac{T_{n+2}(1)-T^2_{n+1}(1)}{m^2T_n(1)}-\frac1m\]
 to derive the asymptotics  
 \[{\bf E} X_n \sim\frac n{\log n},\quad {\bf var}(X_n)\sim\frac n{\log^2n}.\]
 The last step relied on the asymptotic analysis of $T_n(1)$. While the calculations for the expected value were a straightforward application of the saddle point method, the variance was more delicate due to cancellations in $T_{n+2}(1)-T^2_{n+1}(1)$.  
 
 Alternatively, the asymptotic normality follows from  Theorem~\ref{thm} with $\g(x)=x+c$ and $c(x)=0$.
 
 \subsection{ Whitney numbers, Stirling numbers and their generalizations}
 The coefficients $(T_{n,k})$ of the polynomials  $(T_n(x))$ given by \eqref{wang_poly}
 satisfy the recurrence
 \be\label{w-rec}T_{n,k}=T_{n-1,k-1}+(mk+c)T_{n-1,k},\quad 0\le k\le n.\ee
 Versions of numbers satisfying \eqref{w-rec} frequently appear in the literature under various names. In particular, when $c=1$ they are referred to as Whitney numbers of the second kind \cite{Be1, Be2}. The (ordinary) generating functions of Whitney numbers are called Dowling polynomials and the row sums of Whitney numbers are known as Dowling numbers. When $m=2$ Whitney numbers appear as A039755 (and A039756) in \cite{oeis}  under the name B--analogs of  Stirling numbers of the second kind.  Sequences A007405, A003575--A003582, A364069 and A364070 are Dowling numbers for $m=2,3,\dots, 10$,  $m=64$ and $m=624$, respectively.
 
The case $c=0$ are translated Whitney numbers (see, e.g. \cite{bb}).  Examples in \cite{oeis} are sequences A075497 through A075505  which correspond to   $m=2,\dots,10$.  The case $c=r$ are the  $r$--Whitney numbers \cite{cj}. The latter are also referred to as the  $(r,\beta)$--Stirling numbers \cite{cca} ($c=r$, $m=\beta$). This is because the numbers for the case $\beta=1$ are essentially the $r$--Stirling numbers of the second kind  \cite{br}. Specifically, the $r$--Stirling number $\left\{n\atop k\right\}_r$ counts the number of partitions of the set $[n]:=\{1,2,\dots,n\}$  into $k$ blocks, such that the numbers $1$ through  $r$ are in different blocks. If $T_{n,k}$ are $(r,1)$--Stirling numbers as defined in \cite{cca} then the relation is 
 \[\left\{n\atop k\right\}_r=T_{n-r,k-r} \quad \mbox{for}\quad n\ge k\ge r.\]
 Since $T_{n,k}$  satisfy \eqref{w-rec} with $m=1$,  $c=r$, all $r$--Stirling numbers $\left\{n\atop k\right\}_r$ satisfy \eqref{w-rec} with $m=1$, $c=0$
 but with different initial condition, namely, $\left\{n\atop k\right\}_r=\delta_{n,r}$, $n\le r$. 
 For $r=2, 3, 4$, $r$--Stirling numbers are sequences A143494--A143496 in \cite{oeis}.  Of course, classical Stirling numbers of the second kind correspond to $r=1$ and their $(n/\log n,n/\log^2n)$--asymptotic normality is well--known and was established by Harper \cite{h}, see also a discussion in \cite[Example~III.11 and Proposition~IX.20]{fs}[. By Theorem~\ref{thm} all the variants mentioned above follow the same distribution.
 
Whitney numbers were introduced in the context of geometric lattices associated with groups  \cite{d}, see also \cite{Be0}. Later, combinatorial interpretations (mostly related to restricted and colored set partitions) were found. A general nature of these restrictions is discussed in  \cite{gn}.  There seem to be overlaps in the literature between these various families. Partially for this reason,  we limited references to the papers most relevant here. More details and history may be found therein.

\subsection{Further examples} 
 Other examples of sequences in \cite{oeis} satisfying \eqref{w-rec} are A186695 ($m=2$, $c=-1$) or A111577 ($m=3$, $c=-2$). Both are referred to as Galton triangles and both have $c<0$. 
 
 Sequences satisfying \eqref{w-rec} with $c=m-1$ are referred to as (scaled) Stirling--Frobenius subset numbers. For $m=1$ through $m=4$ they are  A048993, A039755, A225468 and A225469 in \cite{oeis}, respectively. 
 
 Numbers $(S2[d,a](n,k))$ where $a$, $d$ are non--negative integers with $\mbox{gcd}(d,a)=1$ are called Sheffer triangles (see \cite{L2} for a general discussion and \cite[Example~4]{L1} for an example relevant here). They satisfy the recurrence
 \[S2[d,a](n,k)=dS2[d,a](n-1,k-1)+(a+dk)S2[d,a](n-1,k).
  \]
  Thus, their generating polynomials satisfy  \eqref{rec} with $\g(x)=dx+a$, $m=d$ and $c(x)=0$. Sequences  A048993, A039755, A154537, A282629, A225466, A285061 and A225467--A225469 are the numbers $S2[d,a]$ for various values of $d$ and $a$.  
  The scaled Stirling--Frobenius subset numbers mentioned earlier are special cases $S2[m,m-1]$.  We note that the \lq\lq unscaled\rq\rq Stirling--Frobenius numbers are numbers satisfying \eqref{w-rec} with $c=m-1$. This situation was discussed earlier.
    
  While some of these families of numbers appeared in different contexts, from the point of view of the asymptotics, their behavior is the same. By Theorem~\ref{thm} they are all asymptotically normal with the mean $n/\log n$ and the variance $n/\log^2n$. For some of these numbers their asymptotic normality has been explicitly stated;  for others it has not, it seems. However, it should be noted  that the bivariate generating function is available in the explicit form (and has been derived in some cases).  This could then be used to carry out the asymptotic analysis.   The methods for  deriving the bivariate generating function varied from case to case. The method of characteristics, outlined above, seems to provide some uniformity in this respect.

 \subsection{Set partitions without small blocks: $s$--associated Stirling numbers}
 In most of the cases in the literature $c(x)$ is identically zero. One natural example where this is not the case is provided by 
  the recurrence for $0\le k\le n$:
\[D_{n,k}=kD_{n-1,k}+(n-1)D_{n-2,k-1}, \quad D_{0,0}=1.
\] 
  The  numbers $D_{n,k}$ count the number of set partitions of an $n$--set into $k$ blocks, each of size at least 2  (see e.g. \cite{bm} for a relatively recent reference and discussion of some of its properties).   It results in the following recurrence for the generating polynomials
\[
D_n(x)=xD_{n-1}'(x)+(n-1)xD_{n-2}(x),\quad D_0(x)=1.
\]
Theorem~\ref{thm} applies with $\g(x)=0$, $m=1$ and $c(x)=x$.   Thus, for partitions of an $n$--set into blocks of size at least two,  the number of blocks is asymptotically normal with the expected value asymptotic to $n/\log n$ and the variance asymptotic to $n/\log^2n$. In \cite{bm}, real--rootedness of the polynomials $D_n(x)$ was established and its various consequences have been discussed although the asymptotic normality was not one of them.

This example is actually a special case of the so--called $s$--associated Stirling numbers of the second kind (see \cite[pp.~221--222]{C}). They count the number of set partitions into blocks of sizes at least $s$.   The analogous recurrence is 
\[D_{n,k}=kD_{n-1,k}+\binom{n-1}{s-1}D_{n-s,k-1}.\]
This gives the polynomial recurrence
\[D_n(x)=xD'_{n-1}(x)+\frac x{(s-1)!}(n-1)_{s-1}D_{n-s}(x),\]
where $(x)_m$ is the falling factorial. When $s\ge3$ this is technically outside of the scope of Theorem~\ref{thm}, but one can argue in exactly the same way: for 
\[F(z,x)=\sum_n\frac{z^n}{n!}D_n(x)\]
  PDE \eqref{pde} takes the form:
 \[\frac\partial{\partial z}
 F(z,x)-x\frac\partial{\partial x}F(z,x)=\frac x{(s-1)!}z^{s-1}F(z,x),\]
  and hence, with $x(z)=\xi e^{-z}$, 
 \[F(z,x(z))=\exp\left\{\frac\xi{(s-1)!}\int z^{s-1}e^{-z}dz\right\}=\exp\Big\{\xi\Big(-\sum_{j=1}^{s}\frac{z^{s-j}}{(s-j)!}e^{-z}+C\Big)\Big\}.
 \]
 As $1=F(0,x)=e^{\xi(-1+C)}$, $C=1$, and eliminating $\xi=xe^z$  gives
 \be\label{bm}F(z,x)=\exp\Big\{x\Big(e^z-\sum_{j=0}^{s-1}\frac{z^j}{j!}\Big)\Big\}\ee
 as given e.g. in \cite{C}.  The asymptotic analysis applies with
 \[Q_1(z,x)=-x\sum_{j=0}^{s-1}\frac{z^j}{j!}\quad\mbox{and}\quad Q_2(z)=z\]
 in \eqref{f-rep}.
 This yields, as for the cases $s=1$ and $s=2$, the asymptotically normal law with mean $n/\log n$ and the variance $n/\log^2n$.
 
 \subsection{Associated $r$--Whitney  numbers} We close with one more example in the same spirit as the previous one. As we mentioned earlier, a recent paper \cite{gn} combined two different restrictions imposed on  set partitions. One  concerns the sizes of parts (association). The other insists that specific elements are in different blocks of a partition. As we will see, our results apply to the number of blocks in such partitions. As far as we know, the number of blocks in such partitions has not been considered before. 

 Following \cite{cj,cca} we say that a set partition is a Whitney colored  $r$--partition 
with $m$ colors if it is a partition of $[r+n]$ such that:
 \begin{itemize}
 \item[(i)]{} the numbers $1,\dots, r$ are in different blocks, 
\item[(ii)]{} the smallest elements of the blocks are not colored,
\item[(iii)]{} elements in  blocks containing $1,\dots, r$ are not colored,
\item[(iv)]{} the remaining elements are colored with $m$ colors.
 \end{itemize}
 Elements $1,\dots,r$ are called distinguished and the blocks containing them are called distinguished blocks. All other blocks and elements are called non--distinguished. In \cite{gn}, the $s$--associated $r$--Dowling numbers $(D_{n,m,r}^{\ge s})$  are defined and some of their properties are studied. Combinatorially,  $D_{n,m,r}^{\ge s}$ is the number of Whitney colored $r$--partitions with $m$ colors with the property that each non--distinguished block contains at least $s$ elements.  In analogy with Dowling numbers being the row sums of Whitney numbers, we let the $s$--associated $r$--Whitney number $W_{n,k,m,r}^{\ge s}$ be the number of such partitions with $k$ non--distinguished blocks. 
 We will show that these numbers are asymptotically normal with the mean $n/\log n$ and the variance $n/\log^2n$. As $m$, $r$ and $s$ are fixed we will write $W_{n,k}=W_{n,k,m,r}^{\ge s}$ through the rest of this section. The numbers $W_{n,k}$ satisfy the recurrence (see \cite[Proof of Theorem~3]{gn} for an argument for the Dowling counterparts)
 \[W_{n,k}=(r+mk)W_{n-1,k}+\binom{n-1}{s-1}m^{s-1}W_{n-s,k-1}.\]
 Indeed, the first term counts the instances in which  $n+r$ is in one of the distinguished blocks or in a non-distinguished block of size larger than $s$. The second counts instances in which it is in a non--distinguished block of size $s$. (In the latter case one needs to chose any $s-1$ elements from $\{r+1,\dots,r+n-1\}$ for that block and  color all but the smallest one in $m^{s-1}$ ways).  Row generating polynomials 
 \[ W_n(x):=\sum_{k=0}^nW_{n,k}x^k\]
 satisfy
 \[W_n(x)=rW_{n-1}(x)+mxW_{n-1}'(x)+m^{s-1}\binom{n-1}{s-1}xW_{n-s}(x).\]
 The resulting PDE for the bivariate generating function $F(z,x)$ takes a slightly different form than in the previous example. Namely,
 \[\frac\partial{\partial z}F(z,x)-mx\frac\partial{\partial x}F(z,x)=\left(r+\frac x{(s-1)!}(mz)^{s-1}\right)F(z,x).\]
 Following the same steps as in that example gives
 \[F(z,x)=\exp\Big\{rz+C_0+\frac\xi m\Big(-\sum_{j=0}^{s-1}\frac{(mz)^{j}}{j!}e^{-mz}+C_1\Big)\Big\},\]
 where $\xi=xe^{mz}$ and $C_0$, $C_1$ are integration constants. The initial condition $F(0,x)=1$ leads to 
  \[F(z,x)=\exp\Big\{rz+\frac xm\Big(e^{mz}-\sum_{j=0}^{s-1}\frac{(mz)^{j}}{j!}\Big)\Big\}.\]
  The aforementioned  asymptotic normality of $(W_{n,k})$ follows by the same asymptotic analysis as before with 
  \[Q_1(x,z)=rz-\frac xm\sum_{j=0}^{s-1}\frac{m^{j}}{j!}z^j\quad\mbox{and} \quad Q_2(z)=\frac zm\] 
  in \eqref{f-rep}.
  
  We note that the univariate exponential generating function of the sequence $(D_{n,m,r}^{\ge s})$ is
  \[F(z,1)=\exp\Big\{rz+\frac 1m\Big(e^{mz}-\sum_{j=0}^{s-1}\frac{(mz)^{j}}{j!}\Big)\Big\},\]
   as was derived in \cite[Theorem~2]{gn} by a different method.
   
   \subsection{Final comments} Recurrences of type \eqref{gen-rec} are very common in combinatorial probability. In vast majority of cases $c_n(x)=0$ and $b_n(x)=b(x)$ so the recurrence simplifies to
   \[P_n(x)=a_n(x)P_{n-1}(x)+b(x)P'_{n-1}(x).\]
   Paper \cite{hcd} comprehensively treated the case $a_n(x)=\a(x)n+\g(x)$ with $\a(1)>0$ and $b_n(x)=(1-x)\beta(x)$. The present work covers  the situation $a_n(x)=\g(x)$ and $b(x)=mx$, $m>0$. However, there are examples of recurrences of the above type with the coefficients $a_n(x)$ and $b(x)$ of different form that those just mentioned. Thus, it would seem worthwhile  to study such recurrences for other sequences of interest. For some cases it should be rather straightforward, for other might be more challenging. In particular, for the method of moments to work well it is very important that  the condition $b(1)=0$ holds.  The approach based on the  characteristics is less sensitive to that requirement. However, its drawback is that the resulting PDE might not have a closed form solution. 
  In such cases one would have to work with implicitly defined function $F(z,x)$ (see \cite[Sections~3.1 and 5.1--5.3]{hcd} for some discussion of that aspect) or develop other approaches to handle such cases. 
   
 Another possible direction for research is to consider more general forms for the terms  $c_n(x)P_{n-2}(x)$. (In fact, the last two examples are of that type.) It is unclear, however,  how common such recurrences are.


\begin{thebibliography}{14}


\bibitem{bb} H. Belbachir and I. E. Bousbaa.  
\newblock Translated Whitney and  $r$--Whitney numbers: a combinatorial approach.
\textit{J. Integer Seq.} 16 (2013), no. 8, Article 13.8.6, 7 pp.

\bibitem{Be0} M. Benoumhani. On Whitney numbers of Dowling lattices. \textit{Discrete Math.} 159: 3--33, 1996. DOI: 10.1016/0012-365X(95)00095-E.

\bibitem{Be1} M. Benoumhani. 
\newblock On some numbers related to Whitney numbers of Dowling lattices. 
\newblock\textit{Adv. in
Appl. Math.}  19(1): 106--116, 1997.
\newblock DOI: 10.1006/aama.1997.0529.
%
\bibitem{Be2} M. Benoumhani. 
\newblock Log-concavity of Whitney numbers of Dowling lattices. 
\newblock\textit{Adv. in Appl. Math.}  22(2):186--189, 1999.
\newblock DOI: 10.1006/aama.1998.0621.
%

\bibitem{bm}
M.~B\'ona and I. Mez\H{o}. 
\newblock Real zeros and partitions without singleton blocks.
\newblock \textit{European J. Combin.} 51: 500--510,  2016.
\newblock DOI: 10.1016/j.ejc.2015.07.021.

\bibitem{br}
A. Z. Broder. The $r$--Stirling Numbers. \textit{Discrete Math.}  49: 
241-259, 1984. DOI: 10.1016/0012-365X(84)90161-4.

\bibitem{cj}
G.--S. Cheon and J.--H. Jung.
$r$--Whitney numbers of Dowling lattices.
\textit{Discrete Math.} 312: 2337–2348,  2012. DOI: 10.1016/j.disc.2012.04.001.

\bibitem{C}
L.~Comtet.
\newblock Advanced Combinatorics.
\newblock Reidel, 1974.
\newblock ISBN: 90-277-0441-4.

\bibitem{cca}
Corcino, R.B., Corcino, C.B., Aldema, R. Asymptotic normality of the $(r,\beta)$--Stirling
numbers. \textit{Ars Combin.} 81:  81–96, 2006.


\bibitem{DDJ}
D.~Dominici, K.~Driver, and K~Jordaan.
\newblock Polynomial solutions of differential--difference equations.
\newblock \textit{J. Approx. Theory} 163: 41-48, 2011.
\newblock DOI: 10.1016/j.jat.2009.05.010.

\bibitem{d}
T.~A.~Dowling. A class of geometric lattices based on finite groups. \textit{J. Combin.Theory, Ser. B} 14:  
61-86, 1973. Erratum: \textit{J.~Combin. Theory, Ser. B 15}, 211, 1973. 

\bibitem{fs}
P.~Flajolet and R.~Sedgewick.
\newblock Analytic Combinatorics.
\newblock Cambridge University Press, 2009.
\newblock ISBN: 978-0-521-89806-5.
  \newblock DOI: 10.1017/CBO9780511801655.

\bibitem{gn}
E. Gyimesi and G. Nyul.
Associated $r$--Dowling numbers and some relatives.
\textit{C.~R.~Math. Acad. Sci. Paris}  359:  47-55, 2021. DOI:10.5802/crmath.145. 

\bibitem{h} L. H. Harper.
     \newblock {S}tirling behaviour is asymptotically normal.
   \newblock\textit{Ann. Math. Statist.} 
  38: 410--414, 1967.
  DOI: 10.1214/aoms/1177698956.
 

\bibitem{hay} W. K. Hayman.
     \newblock A generalisation of {S}tirling's formula.
   \newblock\textit{J. Reine Angew. Math.} 
  196: 67--95, 1956.
  \newblock DOI: 10.1515/crll.1956.196.67.
  
\bibitem{HL}
P.~Hitczenko  and A.~Lohss.
\newblock Probabilistic consequences of some polynomial recurrences. 
\newblock\textit{Random Struct. Alg.}   53: 652--666, 2018.
\newblock DOI: 10.1002/rsa.20820.


\bibitem{hcd}
H.~K.~Hwang.
H.~H.~Chern and G.-H.~Duh.
\newblock An asymptotic distribution theory for Eulerian recurrences with
   applications.
\newblock \textit{Adv.~in~Appl.~Math.} 112 (2020) 101960, 125.
\newblock DOI: 10.1016/j.aam.2019.101960.


\bibitem{pde} F. John.
     \newblock Partial Differential Equations.
    \newblock Springer, 4th edition, 1982.
    \newblock ISBN: 0-387-90609-6.
  \newblock
       DOI: 10.1007/978-1-4684-9333-7.
       
       \bibitem{L1}
       W.~Lang. On generating functions of diagonals sequences of Sheffer and Riordan number triangles. \texttt{arXiv:1708.01421},  2017.
   
   
\bibitem{L2}W.~Lang. 
On sums of powers of arithmetic progressions, and generalized Stirling, Eulerian and Bernoulli numbers. \texttt{arXiv:1707.04451}, 2017.    
       

\bibitem{LW}
L.~L.~Liu and Y.~Wang.
\newblock A unified approach to polynomial sequences with only real zeros.
\newblock \textit{Adv. Appl. Math.} 38: 542--560, 2007.
\newblock DOI: 10.1016/j.aam.2006.02.003.

\bibitem{m} I.~Mez\H{o}.
     \newblock The {$r$}-{B}ell numbers.
   \newblock\textit{J. Integer Seq.} 
  14, Article 11.1.1, 14, 2011.
  
  \bibitem{n} E.~Neuwirth.
\newblock Recursively defined combinatorial functions: Extending Galton's board.
\newblock\textit{ 
Discrete Math.} 239: 33-51, 2001.
\newblock DOI: 10.1016/S0012-365X(00)00373-3.

\bibitem{oeis} OEIS Foundation Inc. (2024). The On-Line Encyclopedia of Integer Sequences. Published electronically at https://oeis.org.

 
 \bibitem{s} R.~Suter.
     \newblock Two analogues of a classical sequence.
   \newblock\textit{J. Integer Seq.}, 3: 2000, Article 00.1.8.

  
\bibitem{wang}
D.~G.~L.~Wang.
\newblock On colored set partitions of type $B_n$.
\newblock \textit{Cent. Europ. J. Math.} 12(9): 1372--1381, 2014.
\newblock DOI: 10.2478/s11533-014-0419-9.



\end{thebibliography}
\end{document}